\documentclass[12pt]{article}                
\usepackage{graphicx}
\usepackage{psfrag}
\usepackage{amsmath, amssymb}
\usepackage{mathptmx}
\usepackage{ntheorem}
\usepackage{mathrsfs}

\usepackage{setspace}
\makeatletter
\newsavebox\myboxA
\newsavebox\myboxB
\newlength\mylenA
\newcommand*\xoverline[2][0.75]{%
    \sbox{\myboxA}{$\m@th#2$}%
    \setbox\myboxB\null
    \ht\myboxB=\ht\myboxA%
    \dp\myboxB=\dp\myboxA%
    \wd\myboxB=#1\wd\myboxA
    \sbox\myboxB{$\m@th\overline{\copy\myboxB}$}
    \setlength\mylenA{\the\wd\myboxA}
    \addtolength\mylenA{-\the\wd\myboxB}%
    \ifdim\wd\myboxB<\wd\myboxA%
       \rlap{\hskip 0.5\mylenA\usebox\myboxB}{\usebox\myboxA}%
    \else
        \hskip -0.5\mylenA\rlap{\usebox\myboxA}{\hskip 0.5\mylenA\usebox\myboxB}%
    \fi}
\makeatother

\DeclareMathOperator{\I}{I}

\newcommand{\e}{\mathbb{E}}

\newcommand{\Reals}{\mathbb{R}}

\newcommand{\la}{\langle}
\newcommand{\ra}{\rangle}
\newcommand\qed{\hfill\hbox{\rlap{$\sqcap$}$\sqcup$}}

\def\Bla{\Big{\langle}}
\def\Bra{\Big{\rangle}}

\newtheorem{theorem}{\bf Theorem}

\newtheorem{remark}{\bf Remark}

\newtheorem{proposition}{\bf Proposition}

\newenvironment{Proof of lemma}{\noindent{\bf Proof of Lemma}}{\hfill$\Box$\newline}
\newenvironment{Proof of theorem}{\noindent{\bf Proof of Theorem}}{\hfill{\footnotesize${\square}$}\newline}
\newenvironment{Proof of theorems}{\noindent{\bf Proof of Theorems}}{\hfill$\Box$\newline}
\newenvironment{Proof of proposition}{\noindent{\bf Proof of Proposition}}{\hfill$\Box$\newline}
\newenvironment{Proof of propositions}{\noindent{\bf Proof of Propositions}}{\hfill$\Box$\newline}
\newenvironment{Proof of exercise}{\noindent{\it Proof of Exercise:}}{\hfill$\Box$}
\theoremstyle{nonumberplain}
\newcommand\specialref{}

\begin{document}

	\nocite{*} 
	
	\title{Temperature chaos in some spherical\\ mixed $p$-spin models}
	
	\author{
		Wei-Kuo Chen\thanks{\textsc{\tiny School of Mathematics, University of Minnesota, wkchen@umn.edu. Partially supported by NSF grant DMS-1642207 and Hong Kong Research Grants Council
				GRF-14302515.}}
		\and
		Dmitry Panchenko\thanks{\textsc{\tiny Department of Mathematics, University of Toronto, panchenk@math.toronto.edu. Partially supported by NSERC.}}
	}
	\date{}
	\maketitle
	
\begin{abstract}
We give two types of examples of the spherical mixed even-$p$-spin models for which chaos in temperature holds. These complement some known results for the spherical pure $p$-spin models and for models with Ising spins. For example, in contrast to a recent result of Subag who showed absence of chaos in temperature in the spherical pure $p$-spin models for $p\geq 3$, we show that even a smaller order perturbation induces temperature chaos.
\end{abstract} 
\vspace{0.5cm}

\section{Introduction}

In a recent paper \cite{Sub16}, building upon earlier work in \cite{AufBen, AufBenCer, Sub15, SubZ15}, Subag proved that there is no chaos in temperature in the spherical pure $p$-spin models for $p\geq 3$, at low enough temperature; this result was obtained as a consequence of a detailed geometric-probabilistic description of the support of the Gibbs measure in these models. Spherical pure $p$-spin models are believed to be one of a few special cases for which chaos in temperature does not hold (another example is in Proposition \ref{prop2} below), and one expects chaos in temperature for many spherical mixed $p$-spin models, as well as for models with Ising spins. In the case when the mixture is one-step replica symmetry breaking, and at low enough temperature, chaos in temperature can be proved by an adaptation of the techniques in \cite{Sub16}; this will appear in the future work, \cite{Sub16b}. In this paper, we will give two types of examples of spherical mixed $p$-spin models for which chaos in temperature holds. The advantage of our results is that they hold at any temperature, and one of the examples is not restricted to one-step replica symmetry breaking case. The disadvantage is that the proofs are purely analytic and do not come with a description of the Gibbs measure beyond what is already known.

For $N\geq 1$, let us denote the sphere of radius $\sqrt{N}$ in $\Reals^N$ by $S_N$ and let $\nu_N$ be the uniform probability measure on $S_N$. For $p\geq 1,$ we consider the spherical pure $p$-spin Hamiltonian
\begin{align}
H_{N,p}(\sigma)&=\frac{1}{N^{(p-1)/2}}\sum_{1\leq i_1,\ldots,i_p\leq N}g_{i_1,\ldots,i_{p}}\sigma_{i_1}\cdots\sigma_{i_{p}},
\label{Hamp}
\end{align}
where $\sigma\in S_N$ and $(g_{i_1,\ldots,i_{p}})$ are i.i.d. standard Gaussian random variables for all $i_1,\ldots,i_{p}$ and $p\geq 1$. The Hamiltonian of the mixed $p$-spin model is defined as a linear combination
\begin{equation}
H_N(\sigma)=\sum_{p\geq 1}\gamma_pH_{N,p}(\sigma),
\label{Hammixed}
\end{equation}
where, to ensure that the series is well defined, we assume that $\sum_{p\geq 1}2^p\gamma_p^2<\infty$. The covariance of this Hamiltonian is given by
\begin{equation}
\e H_N(\sigma^1)H_N(\sigma^2)=N\xi\bigl(R(\sigma^1,\sigma^2)\bigr),
\end{equation}
where the function $\xi(x)=\sum_{p\geq 1}\gamma_p^2 x^{p}$ and 
\begin{equation}
R(\sigma^1,\sigma^2)= \frac{1}{N}\sum_{i\leq N} \sigma_i^1\sigma_i^2
\end{equation}
is the overlap of $\sigma^1$ and $\sigma^2.$ From now on, we only consider mixed \emph{even-$p$-spin models}, that is,
\begin{equation}
\gamma_p = 0 \mbox{ for all odd } p\geq 1.
\end{equation}

For a given inverse temperature parameter $\beta>0,$ we recall the definitions of the free energy and the partition function,
\begin{equation}
F_{N,\beta} =\frac{1}{N}\e \log Z_{N,\beta} \mbox{ and }
Z_{N,\beta}=\int_{S_N}\exp\beta H_N(\sigma)\, \nu_N(d\sigma),
\end{equation}
as well as the Gibbs measure,
\begin{align}\label{gibbs}
G_{N,\beta}(d\sigma)&=\frac{\exp\beta H_N(\sigma)}{Z_{N,\beta}}\, \nu_N(d\sigma).
\end{align}
Given two inverse temperature parameters $\beta_1,  \beta_2 >0$, we will denote by $(\tau^\ell,\rho^\ell)_{\ell\geq 1}$ the i.i.d. sample from the product measure $G_{N,\beta_1}\times G_{N,\beta_2}$, and we will use the standard notation $\la\, \cdot\,\ra$ for the Gibbs average with respect to $(G_{N,\beta_1}\times G_{N,\beta_2})^{\otimes\infty}$. In this paper, chaos in temperature means that, for $\beta_1\not = \beta_2$,
\begin{equation}
\lim_{N\rightarrow\infty}\e \bigl\la \bigl|R(\tau^1,\rho^1)\bigr| \bigr\ra=0.
\label{chaosinT}
\end{equation}
This limit essentially says that $\tau^1$ and $\rho^1$ are orthogonal to each other if there is a different between two temperatures. 

We will describe two examples when chaos in temperature holds. In order to state our main results, we first need to recall the definition of the Parisi measure. It was proved in \cite{TalagrandSphh, C13} that the limit of the free energy can be computed through the Crisanti-Sommers formula \cite{CS},
\begin{align}\label{cs}
\lim_{N\rightarrow\infty}F_{N,\beta}&=\inf_{\alpha\in\mathcal{M}}\mathcal{Q}_\beta(\alpha),
\end{align}
where $\mathcal{M}$ is the collection of all cumulative distribution functions $\alpha$ on $[0,1]$ with $\alpha(\hat{s})=1$ for some $\hat{s}<1,$ and the functional $\mathcal{Q}_\beta$ is defined as
\begin{align}
\mathcal{Q}_\beta(\alpha)&=\frac{1}{2}\Bigl(\beta^2\int_0^1\xi'(s)\alpha(s)\, ds+\int_0^{\hat{s}}\frac{ds}{\int_s^1\alpha(q)dq}+\log(1-\hat{s})\Bigr).
\label{cs2}
\end{align}
This functional is well-defined and independent of the choice of $\hat{s}.$ It is also strictly convex and continuous with respect to the $L_1$-distance on $\mathcal{M}$. Note that $\mathcal{M}$ is convex, but is not compact. Nonetheless, it can be shown (see e.g. \cite{TalagrandSphh} or \cite{JT}) that the infimum in \eqref{cs} is uniquely achieved by some $\alpha$, which we will denote by $\alpha_\beta$. The probability measure $\mu_\beta$ with the c.d.f. $\alpha_\beta$ is called the Parisi measure. We say that $\mu_\beta$ is replica symmetric (RS) if its support is one point; $\mu_\beta$ is one-step replica symmetric breaking (1-RSB) if its support consists of two distinct points; $\mu_\beta$ is full replica symmetric breaking (FRSB) if its support is not a finite set.

In our first result, we consider a pure even-$p_0$-spin model with asymptotically vanishing perturbation. Consider any two different even integers
\begin{equation}
p_0\geq 4, \mbox{ and } p\not = p_0
\end{equation}
and let us fix any real number $a$ in the interval
\begin{equation}
0<a<\frac{1}{4}.
\label{parameterA}
\end{equation}
Let us consider the Hamiltonian of the form
\begin{equation}
H_N(\sigma) = {H}_{N,p_0}(\sigma)+\frac{1}{N^{a}}\gamma_p H_{N,p}(\sigma).
\label{Hampert}
\end{equation} 
Because of the factor $N^{-a}$, the second term is of a smaller order and can be viewed as a vanishing perturbation of the pure $p_0$-spin Hamiltonian. As a result, it does not affect the limit of the free energy in (\ref{cs}), so the functional defined in (\ref{cs2}) is expressed in terms of $\xi(x) = x^{p_0},$ from which it is known in \cite{TalagrandSphh} that the Parisi measure is either RS or  1-RSB. In contrast to the result of Subag \cite{Sub16}, we will show that this vanishing perturbation term induces temperature chaos. In particular, this indicates that, at least in this case, temperature chaos or its absence can not be detected by free energy calculations.

\begin{theorem}\label{thm1}
If $\beta_1\not =\beta_2$, then there exists $\gamma_p = \gamma_{N,p}\in [1,2]$ possibly varying with $N$ such that chaos in temperature (\ref{chaosinT}) holds.
\end{theorem}
In this result, the fact that $p_0$ and $p$ are even is needed for technical purposes, as we will be using a key result, Theorem 4 in \cite{PT}, which is known only when $p_0$ is even. See Remark \ref{rmk1} below for further details.

To the pure $p$-spin Hamiltonian in the perturbation term in (\ref{Hampert}), one could also add an arbitrary mixed $p$-spin Hamiltonian containing at least one even $p$-spin interaction, if one so wishes. Let us also mention that the condition $a<1/4$ in (\ref{parameterA}) is a standard technical condition to ensure the validity of the Ghirlanda-Guerra identities \cite{GG} for the $p^{\rm th}$ moment of the overlaps (see, e.g., Section 3.2 in \cite{SKmodel}), and is likely not optimal.

Our second example will be in the setting of the so-called generic models. We will call the mixed even-$p$-spin Hamiltonian (\ref{Hammixed}) \emph{generic} if the linear span of functions $x^p$ for even $p\geq 2$ such that $\gamma_p\not = 0$ and constants is dense in $C([0,1],\|\,\cdot\,\|_\infty)$. Denote the smallest point in the support of the Parisi measure by
\begin{equation}
c_\beta = \inf \mathrm{supp}\, \mu_\beta.
\label{smallestcs}
\end{equation}
Given two inverse temperatures $\beta_1, \beta_2>0$, let
\begin{equation}
q_0(\beta_1,\beta_2) = \inf\bigl\{t \,:\, \beta_1 \mu_{\beta_1}\bigl([0,t)\bigr) \not= \beta_2\mu_{\beta_2}\bigl([0,t)\bigr)\bigr\}.
\label{qnot}
\end{equation}
We will need the following condition on the two temperatures,
\begin{equation}
q_0(\beta_1,\beta_2) \leq \max\bigl(c_{\beta_1},c_{\beta_2}\bigr).
\label{q0uncoupled}
\end{equation}
This means that either $c_{\beta_1}\not = c_{\beta_2}$ or, otherwise, the scaled Parisi measures $\beta_1\mu_{\beta_1}$ and $\beta_2\mu_{\beta_2}$ are immediately different to the right of the smallest point in their support $c_{\beta_1}=c_{\beta_2}.$ If this condition holds then, as in \cite{P15}, we say that the Parisi measures $\mu_{\beta_1}, \mu_{\beta_2}$ are \emph{uncoupled}.  The following is our second main result.

\begin{theorem}
	\label{thm2}
	Suppose that the model is generic. If $\beta_1\not = \beta_2$, the condition (\ref{q0uncoupled}) is satisfied and $\min(c_{\beta_1}, c_{\beta_2}) = 0$, then chaos in temperature (\ref{chaosinT}) holds.
\end{theorem}
The proof of this theorem also works for generic models that include both even and odd $p$-spin interactions, in which case one needs a technical assumption that the linear span of the functions $x^p$ for $p\geq 1$ such that $\gamma_p\not = 0$ and constants is dense in $C([-1,1],\|\,\cdot\,\|_\infty)$. We comment more on this right after the proof of Theorem \ref{thm2}.

In a recent work \cite{JT}, Jagannath and Tobasco showed that the problem of computing the Parisi measure in the spherical models can be reduced to a certain finite dimensional optimization problem, see Corollary 1.5 in \cite{JT}. This means that one should be able to easily check the conditions (\ref{q0uncoupled}) and $\min(c_{\beta_1}, c_{\beta_2}) = 0$ in Theorem \ref{thm2} numerically. Their result (see also Theorem 6 in \cite{AufChen1}) implies that these two conditions are equivalent to the following:
\begin{enumerate}
\item 
either $\min(c_{\beta_1}, c_{\beta_2}) = 0$ and $\max(c_{\beta_1}, c_{\beta_2}) > 0$ or, otherwise,

\item $\beta_1 \mu_{\beta_1}(\{0\}) \not= \beta_2\mu_{\beta_2}(\{0\}).$
\end{enumerate}
In the case when the model is RS or 1-RSB at all temperatures, checking the conditions of Theorem~\ref{thm2} is particularly easy, as will be discussed in the next section. 

It is possible that the condition \eqref{q0uncoupled} is not spurious. In the next section, we will give examples of mixed $p$-spin models whose Parisi measures are FRSB and for which \eqref{q0uncoupled} is violated. Unlike in the pure $p$-spin model, in this case it is very challenging to obtain useful control of the free energy to say anything about the cross overlap $|R(\tau^1,\rho^1)|$, and it has been conjectured in \cite{R02} that there is no temperature chaos.

\section{$1$-RSB and FRSB solutions} \label{Sec2label}

In this section, we will first review several known results about the models whose Parisi measures are 1-RSB. These will be useful to us in the proof of Theorem \ref{thm1}. As one shall see, the problem of checking the conditions in Theorem \ref{thm2} simplifies in this case quite a bit. After that, we will describe a criterion that guarantees that the Parisi measure is FRSB and explain how the inequality \eqref{q0uncoupled} can be violated.

By Proposition 2.2 \cite{TalagrandSphh}, if the function $\xi''(s)^{-1/2}$ is convex then the support of the Parisi measure $\mu_\beta$ contains at most two points. In the case of the pure $2$-spin model with $\xi(x)=x^2,$ the same proof actually shows that the Parisi measure is always concentrated on one point. On the other hand, Proposition 2.3 in \cite{TalagrandSphh} gives that, whenever
\begin{equation}\label{high}
\sup_{0<s<1}\bigl(\beta^2\xi(s)+\log(1-s)+s\bigr)\leq 0,
\end{equation}
the Parisi measure is concentrated at $0$, $\mu_\beta=\delta_0$, and, in the complementary case,
	\begin{align}\label{prop1:eq1}
	\sup_{0<s<1}\bigl(\beta^2\xi(s)+\log(1-s)+s\bigr)>0,
	\end{align}
the Parisi measure is not concentrated at $0$, $\mu_\beta\not =\delta_0$. In this case, unless the model is pure $2$-spin, if the Parisi measure has at most two atoms, then it must be of the form
\begin{align}
\label{eq1}
\mu_\beta=m \delta_0+(1-m)\delta_{q}
\,\,\mbox{ for $0< m<1$ and $q>0.$} 
\end{align}
To summarize, we have the following proposition. As its proof does not seem to appear in the literature, we will present a detailed argument in the last section. 
\begin{proposition}\label{prop1}
Suppose $\gamma_p \not = 0$ for some $p\geq 3$. If (\ref{prop1:eq1}) holds and the Parisi measure has at most two atoms, then it is of the form (\ref{eq1}).
\end{proposition}
If $\alpha_{m,q}(s)=m 1_{[0,q)}(s)+ 1_{[q,1]}(s)$ is the c.d.f. of $\mu_\beta$ in (\ref{eq1}), from the optimality of $\mu_{\beta}$, it is easy to check by a direct differentiation of $\mathcal{Q}_{\beta}(\alpha_{m,q})$ with respect to $m$ and $q$ that
	\begin{align}
	\begin{split}
	\label{optimq}
	\beta^2\xi'(q)&=\frac{1}{m}\Bigl(\frac{1}{1-q}-\frac{1}{1-q+mq}\Bigr),\\
	\beta^2\xi(q)&=\frac{1}{m^2}\log\Bigl(\frac{1-q+q m}{1-q}\Bigr)-\frac{q}{m}\frac{1}{1-q+m q}.
	\end{split}
	\end{align}
To check the conditions in Theorem \ref{thm2}, one needs to show that the parameters $m$ corresponding to two different temperatures $\beta_1\not = \beta_2$ satisfy 
\begin{equation}
\beta_1 m_{1} \not = \beta_2 m_{2}.
\label{mbeta}
\end{equation}
This is always true for pure even-$p$-spin models, $\xi(p)=x^p$, for $p\geq 4$, as can be easily checked. As a result, for any given $\beta_1\not = \beta_2$, a small enough generic perturbation of a pure $p$-spin model, for which $\xi''(s)^{-1/2}$ is convex will still satisfy (\ref{mbeta}). Let us now recall several facts in the setting of the pure even-$p$-spin models for $p\geq 4$, which will be useful to us in the proof of Theorem \ref{thm1}. 

First of all, when (\ref{prop1:eq1}) holds, the optimal parameter $m$ is strictly positive, $0<m<1$, so the Parisi measure has two atoms (see, e.g., Section 4 in \cite{PT}). For these models, it is also well known that the Parisi measure $\mu_\beta$ is the limiting distribution of the overlap $R(\sigma^1,\sigma^2).$ Indeed, Theorem~4 in \cite{PT} (applied to two systems at the same temperature) shows that the limiting distribution of $|R(\sigma^1,\sigma^2)|$ concentrates on two points $\{0,q\}$, where $q$ is the second atom in (\ref{eq1}), while the proof of Theorem 1.2 in \cite{TalagrandSphh}) gives
\begin{equation}
\lim_{N\to\infty}\e\bigl\la R(\sigma^1,\sigma^2)^p\bigr\ra = \int s^p \, \mu_{\beta}(ds).
\end{equation}
Clearly, for $\mu_\beta$ as in (\ref{eq1}), these two facts imply that the distribution of $|R(\sigma^1,\sigma^2)|$ converges to $\mu_\beta.$ At two different temperatures, Theorem 4 in \cite{PT} gives that the limiting distribution of the cross-overlap $|R(\tau^1,\rho^1)|$ concentrates on two points $\{0,\sqrt{q_1 q_2}\}$, where $q_1$ and $q_2$ are the non-zero atoms corresponding to these two temperatures. All these facts are proved by free energy calculations, which are not affected by the perturbation term in the Hamiltonian (\ref{Hampert}). Consequently, they all hold for the perturbed model \eqref{Hampert}. These will be used in the proof of Theorem \ref{thm1}.

While the condition that $\xi''(s)^{-1/2}$ is convex guarantees that the Parisi measure is at most $1$-RSB, the next proposition shows that if $\xi''(s)^{-1/2}$ is concave then the Parisi measure is FRSB.
\begin{proposition}
	\label{prop2}
	Suppose that $\xi''(s)^{-1/2}$ is concave on $(0,1]$. If $\beta\xi''(0)^{1/2}>1,$ then 
	\begin{align}
	\alpha_\beta(t)&=\left\{
	\begin{array}{ll}
	\frac{\xi'''(t)}{2\beta\xi''(t)^{3/2}},&\mbox{if $t\in[0,q)$},\\
	1,&\mbox{if $t\in[q,1]$},
	\end{array}
	\right.
	\label{aFRSB}
	\end{align}
	where $q\in(0,1)$ is the unique solution of
	\begin{align}
	\label{prop2:eq1}
	\frac{1}{\beta\xi''(q)^{1/2}}=1-q.
	\end{align}
\end{proposition}

This proposition is  a general statement of the examples discussed after Proposition 2.2 in \cite{TalagrandSphh}, where the author constructed a generic even $p$-spin model and a variant whose Parisi measure has no jump at the top of the support. For a simple example when the assumptions in Proposition \ref{prop2} hold (see Example 4 in \cite{AufChen1}), take $\xi(t)=(1-c)t^2+ct^p$ for any $c>0$ such that 
$$
\frac{c}{1-c}\leq \frac{4(p-3)}{(p-1)p^2}\,\,\mbox{ and }\,\,\frac{1}{2(1-c)}<\beta^2.
$$
We will see in the proof of Proposition \ref{prop2} that the assumptions on $\xi$ ensure that (\ref{aFRSB}) is a well-defined c.d.f., which has non-zero jump at $t=q$. Let us take two inverse temperatures such that
$$
\xi''(0)^{-1/2}<\beta_1<\beta_2.
$$
Note that since the graphs of $(\beta_1\xi'')^{-1/2}$ and $(\beta_2\xi'')^{-1/2}$ are concave and $(\beta_1\xi'')^{-1/2}>(\beta_2\xi'')^{-1/2}$, the corresponding solutions $q_1,q_2$ of \eqref{prop2:eq1} satisfy $0<q_1<q_2$. The form of the c.d.f. in (\ref{aFRSB}) implies that 
$$
\beta_1\alpha_{\beta_1}(t)=\beta_2\alpha_{\beta_2}(t)
\mbox{ for all $t\in [0,q_1).$}
$$
Moreover, because of the jump discontinuity of $\alpha_{\beta_1}(t)$ at $t=q_1$,
$$
\beta_1\alpha_{\beta_1}(q_1)> \beta_2\alpha_{\beta_2}(q_1).
$$ 
Recalling the definition (\ref{qnot}), this implies that $q_0(\beta_1,\beta_2)=q_1>0.$ If \eqref{Hammixed} is not pure $2$-spin then $\xi'''(t)>0$ for $t>0$ and $c_{\beta_1} =c_{\beta_2} = 0,$ so the condition \eqref{q0uncoupled} is violated.

\section{Proof of main results} \label{Sec3label}

\textbf{Proof of Theorem \ref{thm1}.}
Our approach to proving Theorem \ref{thm1} will be based on the Ghirlanda-Guerra identities for the coupled systems as implemented in \cite{CP, ChenChaos}, as well as the consequences of  free energy calculations for the pure even-$p_0$-spin models for $p_0\geq 4$ in \cite{PT}, mentioned in the previous section.

First, from the Crisanti-Sommers formula \eqref{cs}, it was known in \cite{TalagrandSphh} that the limiting free energy is differentiable in $\beta$. Using this differentiability, one can obtain concentration  of the Hamiltonian $H_{N,p_0}(\sigma)$,
\begin{align}\label{eq:add1}
\e\Bla\Bigl|\frac{H_{N,p_0}(\sigma)}{N}-\e\Bla\frac{H_{N,p_0}(\sigma)}{N}\Bra\Bigr|\Bra= o(1),
\end{align}  
by a standard argument (see, e.g., \cite{PGGmixed}, or Section 4 in \cite{AufChen}). Under the assumption $a>0$ in (\ref{parameterA}), a standard application of the concentration of the free energy (see, e.g., Theorem 3.3 in \cite{SKmodel}) shows that, for some choice of $\gamma_p = \gamma_{N,p}\in [1,2]$ possibly varying with $N$,
\begin{align}\label{eq:add2}
\e\Bla\Bigl| \frac{H_{N,p}(\sigma)}{N}-\e\Bla\frac{H_{N,p}(\sigma)}{N}\Bra\Bigr|\Bra=O\Bigl(\frac{1}{N^{1/4}}\Bigr).
\end{align}
It should be emphasized that while the validity of \eqref{eq:add1} relies on the differentiability of the limiting free energy, the display \eqref{eq:add2} follows a general principle and does not depend on whether the model is spherical or not. Now, following an identical argument as in \cite{CP, ChenChaos} (see e.g. Lemma 2 in \cite{CP}), for a bounded function $f$ of the overlaps $R(\tau^\ell,\tau^{\ell'})$, $R(\rho^\ell,\tau^{\ell'})$ and $R(\tau^\ell,\rho^{\ell'})$ of the first $n$ replicas, $\ell, \ell'\leq n$, if we use the integration by parts to the centered Hamiltonian against $f$ as a test function and apply the concentrations \eqref{eq:add1} and \eqref{eq:add2}, we will get, for $\kappa = \beta_2/\beta_1$ and $\phi(x)=x^{p_0}$ or $\phi(x)=x^p$,
\begin{align}
\begin{split}\label{eq:add3}
&\e \la f \phi(R(\tau^1,\tau^{n+1})) \ra+\kappa\e \la f \phi(R(\tau^1,\rho^{n+1}))\ra\\
&\approx \frac{1}{n}\e \la  f\ra\e \la \phi(R(\tau^1,\tau^2))\ra+\frac{1}{n}\sum_{\ell=2}^n \e \la f \phi(R(\tau^1,\tau^\ell))\ra+\frac{\kappa}{n}\sum_{\ell=1}^n \e \la f \phi(R(\tau^1,\rho^\ell))\ra
\end{split}
\end{align}
and
\begin{align}
\begin{split}\label{eq:add4}
&\e \la f \phi(R(\rho^1,\rho^{n+1}))\ra+\frac{1}{\kappa}\e \la f \phi(R(\rho^1,\tau^{n+1}))\ra\\
&\approx \frac{1}{n}\e \la  f\ra\e \la \phi(R(\rho^1,\rho^2))\ra+\frac{1}{n}\sum_{\ell=2}^n \e \la f \phi(R(\rho^1,\rho^\ell))\ra+\frac{1}{\kappa n}\sum_{\ell=1}^n \e \la f \phi(R(\rho^1,\tau^{\ell}))\ra.
\end{split}
\end{align}
Here, $\approx$ means $o(1)$ for $\phi(x)=x^{p_0}$, and $O(N^{a-1/4})$ for $\phi(x)=x^p$, with an extra factor $N^a$ coming from the factor $N^{-a}$ in front of $H_{N,p}(\sigma)$ in (\ref{Hampert}). By the assumption $a<1/4$ in (\ref{parameterA}), in both cases $\approx$ means $o(1)$ as $N\to\infty$. 

Next, it will be convenient to replace the above approximate identities by their exact analogues in the thermodynamic limit. Let us consider any subsequential limit in distribution of the overlaps
$$
(R(\tau^\ell,\tau^{\ell'}) )_{\ell\neq \ell'\geq 1},\, (R(\rho^\ell,\tau^{\ell'}))_{\ell\neq \ell'\geq 1},\, \mbox{ and } (R(\tau^\ell,\rho^{\ell'}))_{\ell,\ell'\geq 1},
$$
under $\e (G_{N,\beta_1}\times G_{N,\beta_2})^{\otimes\infty}.$ By Theorem 2 in \cite{P15}, there exists a pair $(G_1, G_2)$ of random probability measures on a separable Hilbert space such that this limiting distribution coincides with the distribution under $\e (G_1\times G_2)^{\otimes \infty}$ of the array
\begin{align*}
(\tau^\ell\cdot\tau^{\ell'})_{\ell\neq\ell'\geq 1},\,\,(\rho^\ell\cdot\rho^{\ell'})_{\ell\neq \ell'\geq 1},\,\,(\tau^\ell\cdot\rho^{\ell'})_{\ell,\ell'\geq 1},
\end{align*}
where $(\tau^\ell,\rho^\ell)_{\ell\geq 1}$ is an i.i.d. sample from $G_1\times G_2.$ For simplicity of notation, we will continue to use the notation $\la\,\cdot\,\ra$ also for the average with respect to $(G_1\times G_2)^{\otimes \infty}$. Then the above approximate identities become the following exact identities in the limit,
\begin{align}
\begin{split}
\label{id1}
&\e \la f \phi(\tau^1\cdot\tau^{n+1})\ra+\kappa\e \la f \phi(\tau^1\cdot\rho^{n+1})\ra\\
&=\frac{1}{n}\e \la  f\ra\e \la\phi(\tau^1\cdot\tau^2)\ra+\frac{1}{n}\sum_{\ell=2}^n \e \la f \phi(\tau^1\cdot\tau^\ell)\ra+\frac{\kappa}{n}\sum_{\ell=1}^n \e \la f \phi(\tau^1\cdot\rho^\ell)\ra
\end{split}
\end{align}
and
\begin{align}
\begin{split}\label{id2}
&\e \la f \phi(\rho^1\cdot\rho^{n+1})\ra+\frac{1}{\kappa}\e \la f \phi(\rho^1\cdot\tau^{n+1})\ra\\
&=\frac{1}{n}\e \la  f\ra\e \la \phi(\rho^1\cdot\rho^2)\ra+\frac{1}{n}\sum_{\ell=2}^n \e \la f \phi (\rho^1\cdot\rho^\ell)\ra+\frac{1}{\kappa n}\sum_{\ell=1}^n \e \la f \phi(\rho^1\cdot\tau^{\ell})\ra,
\end{split}
\end{align}
for any bounded function $f$ of $(\rho^\ell\cdot\rho^{\ell'})_{1\leq \ell\neq \ell'\leq n}$, $(\tau^\ell\cdot\tau^{\ell'})_{1\leq \ell\neq \ell'\leq n}$ and $(\tau^\ell\cdot\rho^{\ell'})_{1\leq \ell\neq \ell'\leq n}$ and for $\phi(x)=x^{p_0}$ or $\phi(x)=x^p$. 

If at one of the temperatures $\beta_1$ or $\beta_2$ the Parisi measure $\mu_\beta$ concentrates at zero, a simple application of the Cauchy-Schwarz inequality (see Lemma 2 in \cite{PT}) immediately implies the chaos in temperature. Hence, we will only consider the case when both temperatures satisfy (\ref{prop1:eq1}), and the two Parisi measures are of the form 
$$
\mu_{\beta_j}=m_j \delta_0+(1-m_j)\delta_{q_j}.
$$
As we mentioned in Section \ref{Sec2label}, in this case,
\begin{align}
&
\e \bigl\la \I\bigl(|\tau^1\cdot \tau^2|=0 \,\,\mbox{or}\,\,q_1\bigr)\bigr\ra =1,
\nonumber
\\
&
\e \bigl\la \I\bigl(|\rho^1\cdot \rho^2|=0\,\,\mbox{or}\,\,q_2\bigr)\bigr\ra =1,
\label{thm2:proof:eq3}
\\
&
\e \bigl\la \I\bigl(|\tau^1\cdot \rho^1|=0\,\,\mbox{or}\,\,\sqrt{q_1q_2}\bigr)\bigr\ra =1.
\nonumber
\end{align} 
Next, we will show that \eqref{thm2:proof:eq3}, (\ref{id1}) and (\ref{id2}) imply
	\begin{align}
	\begin{split}
	\label{id1.1}
	\e \la f \phi(\tau^1\cdot\rho^{n+1})\ra&=\frac{1}{n}\sum_{\ell=1}^n \e \la f \phi(\tau^1\cdot\rho^\ell)\ra,
	\end{split}\\
	\begin{split}
	\label{id1.2}
	\e \la f \phi(\tau^1\cdot \tau^{n+1})\ra&=\frac{1}{n}\e \la  f\ra \e \la \phi(\tau^1\cdot\tau^2)\ra+\frac{1}{n}\sum_{\ell=2}^n \e \la f \phi (\tau^1\cdot\tau^\ell)\ra,
	\end{split}
	\end{align}
	and
	\begin{align}
	\begin{split}
	\label{id2.1}
	\e \la f \phi(\tau^{n+1}\cdot\rho^1)\ra&=\frac{1}{n}\sum_{\ell=1}^n \e \la f \phi(\tau^n\cdot\rho^\ell)\ra,
	\end{split}\\
	\begin{split}
	\label{id2.2}
	\e \la f \phi(\rho^1\cdot\rho^{n+1})\ra&=\frac{1}{n}\e \la  f\ra\e \la \phi(\rho^1\cdot\rho^2)\ra+\frac{1}{n}\sum_{\ell=2}^n \e \la f \phi (\rho^1\cdot\rho^\ell)\ra,
	\end{split}
	\end{align}
for any even function $\phi$ on $[-1,1].$ Once we have these identities, the proof of the temperature chaos is identical to the proof of the first case of Theorem 3 in \cite{CP} (in other words, the proof is entirely in terms of the overlaps and does not distinguish between Ising or spherical spins). 
	
    The verification of the above identities runs as follows. Since $|\tau^1\cdot\tau^2|$ is supported by $0$ and $q_1$ and $|\tau^1\cdot\rho^1|$ is supported by $0$ and $\sqrt{q_1q_2}$, we can rewrite \eqref{id1} with $\phi(x)=x^{d}$ for $d=p_0$ or $d=p$ as
    \begin{align}
    &\e \la f \I(|\tau^1\cdot\tau^{n+1}|=q_1)\ra+\kappa\Bigl(\frac{q_2}{q_1}\Bigr)^{d/2}\e \la f \I(|\tau^1\cdot\rho^{n+1}|=\sqrt{q_1q_2})\ra \nonumber\\
    &=\frac{1}{n}\e \la  f\ra\e \la \I(|\tau^1\cdot\tau^2|=q_1)\ra+\frac{1}{n}\sum_{\ell=2}^n \e \la f \I (|\tau^1\cdot\tau^\ell|=q_1)\ra\nonumber\\
    &+\frac{\kappa}{n}\Bigl(\frac{q_2}{q_1}\Bigr)^{d/2}\sum_{\ell=1}^n \e \la f \I(|\tau^1\cdot\rho^\ell|=\sqrt{q_1q_2})\ra.
    \label{almostthere}
    \end{align} 
For $\beta_1\neq \beta_2$ it can be seen from (\ref{optimq}) that $q_1\neq q_2$. Indeed, if we denote $x = mq/(1-q)$ then the ratio of the two equations in (\ref{optimq}) for $\xi(q)=q^p$ can be rewritten as
$$
\frac{1}{p} = \frac{1+x}{x^2}\log(1+x) -\frac{1}{x}.
$$
The right hand side is convex and decreasing  from $1/2$ to $0$ as $x$ varies from $0$ to infinity, so there exists a unique solution $x$. This implies that if $q_1 = q_2$ then $m_1=m_2$, which contradicts (\ref{optimq}) when $\beta_1\neq \beta_2$. Since $p_0\neq p$ and $q_1\neq q_2$, the equation (\ref{almostthere}) can hold simultaneously for $d=p_0$ and $d=p$ only if 
    \begin{align}\label{add:eq1}
    \e \la f \I(|\tau^1\cdot\rho^{n+1}|=\sqrt{q_1q_2})\ra &=\frac{1}{n}\sum_{\ell=1}^n \e \la f \I(|\tau^1\cdot\rho^\ell|=\sqrt{q_1q_2})\ra
    \end{align}
    and
    \begin{align}\label{add:eq2}
    \e \la f \I(|\tau^1\cdot\tau^{n+1}|=q_1)\ra&=\frac{1}{n}\e \la  f\ra\e \la \I(|\tau^1\cdot\tau^2|=q_1)\ra+\frac{1}{n}\sum_{\ell=2}^n \e \la f \I (|\tau^1\cdot\tau^\ell|=q_1)\ra.
    \end{align}
    Consequently, \eqref{id1.1} follows by using the first equation and the fact again that $|\tau^1\cdot\rho^\ell|$ is supported on $\{0,\sqrt{q_1q_2}\},$ and \eqref{id1.2} follows using the second equation and the fact that $|\tau^1\cdot\tau^\ell|$ is supported by $\{0,q_1\}$. The same argument yields \eqref{id2.1} and \eqref{id2.2}.
\qed

\begin{remark}\label{rmk1}
	\rm There are two technical reasons why we need the evenness of $p_0$ and $p.$ First, the equations in \eqref{thm2:proof:eq3} are deduced from the free energy calculations in Theorem 4 in \cite{PT}, where the  evenness of $p_0$ is needed in order to assure the validity of Guerra's replica symmetry breaking bound for the coupled free energy. It remains an open question whether Theorem 4 in \cite{PT} also holds with odd $p_0.$ Second, the key equation \eqref{almostthere} holds for even $p_0$ and $p.$ If one replaces $p_0$ or $p$ by odd integers, \eqref{almostthere} will become different identities, for which it is unclear to the authors how to employ the same argument to obtain \eqref{add:eq1} and \eqref{add:eq2}.
\end{remark}

\medskip
\noindent
\textbf{Proof of Theorem \ref{thm2}.} The proof of the main result in \cite{P15} in the setting of generic models with Ising spins did not distinguish between models with Ising spins or spherical spins up to and including Section 6. All results up to that point were based entirely on the Ghirlanda-Guerra identities for the overlaps in the form (\ref{id1}) and (\ref{id2}). The main idea was to represent these identities as a certain two-system analogue of the invariance principle developed in the proof of the Parisi ultrametricity conjecture in \cite{PUltra} and, as a consequence, deduce crucial information about the overlaps via various duplication and joint clustering properties. Therefore, one only needs to observe that Theorem 14 in Section 6 of \cite{P15} there implies Theorem \ref{thm2}. 

The only comment that must be made is about Theorem 11 in \cite{P15}. Its proof appealed to one result in the setting of the models with Ising spins, which is not available in the spherical models, namely, Theorem 4 in Chen \cite{ChenChaos}; this was done to reduce to the case with external field, when the overlap of two configurations at the same temperature is asymptotically nonnegative by Theorem 14.12.1 in Talagrand \cite{SG2-2}. However, this was for convenience only and is not necessary. Using the spin flip symmetry of even-$p$-spin models without external field, one can modify the statement of Theorem 11 in \cite{P15} to $\e\la\I(|\sigma^1\cdot\rho^1| = |\sigma^1\cdot\rho^2|)\ra =1$, i.e., for the absolute values of the overlaps, with no changes in the proof (one can appeal to spin flip symmetry to ensure that we can choose the parameter $x$ there to be positive). Theorem 11 was used only in this form (for the absolute values) in all results that follow, so the proof of Theorem 14 is unchanged.
\qed

\medskip
\noindent

\begin{remark}
	\label{rmk2}\rm
	For generic models that include odd-spin interactions we do not have spin flip symmetry, but, instead, we have full set of the Ghirlanda-Guerra identities and, as a consequence, the Talagrand positivity principle (see \cite{SG}, \cite{SG2-2}, or \cite{SKmodel}). This again ensures that the parameter $x$ in the proof of Theorem 11 in \cite{P15} is positive and no modifications are necessary.
\end{remark}

\medskip
\noindent
\textbf{Proof of Proposition \ref{prop1}.} Suppose that the support of $\mu_\beta$ contains at most two points, say $u,v$ with $0\leq u\leq v<1$. We claim that $u=0$. Assume that $u>0.$ For any $\alpha\in\mathcal{M}$ with $\alpha(v)=1$, the optimality of $\alpha_\beta$ gives, by a direct computation (see, e.g., Lemma 2.1 in \cite{TalagrandSphh}) of the derivative of $\mathcal{Q}_\beta$ along the linear path from $\alpha$ to $\alpha_\beta$, 
	\begin{align}\label{prop1:proof:eq1}
	\int_0^{v}\bigl(\alpha(s)-\alpha_\beta(s)\bigr)\Bigl(\beta^2\xi'(s)-\int_0^s\frac{dt}{\bigl(\int_t^1\alpha_\beta(r)dr\bigr)^2}\Bigr)ds\geq 0.
	\end{align}
	For $t\in (0,v)$, define $\alpha_t\in \mathcal{M}$ by
	$$
	\alpha_t(s)=\alpha_\beta(u)1_{[t,v)}(s)+\alpha_\beta(v)1_{[v,1]}(s).
	$$
	Note that the second bracket in the integral of \eqref{prop1:proof:eq1} is continuous. Plugging $\alpha_t$ into \eqref{prop1:proof:eq1} and sending $t\rightarrow u^\pm$ yield that
	\begin{align}\label{Cu}
	\beta^2\xi'(u)=\int_0^u\frac{dt}{\bigl(\int_t^1\alpha_\beta(r)dr\bigr)^2}=Cu,
	\end{align}
	for $C:=\bigl(\int_0^1\alpha_\beta(r)dr\bigr)^{-2}$, where the second equality used the fact that $\alpha_\beta=0$ on $[0,u].$ 
Using (\ref{Cu}), we can write
	\begin{align*}
	\beta^2\xi'(s)-Cs&=\beta^2s\Bigl(\frac{\xi'(s)}{s}-\frac{\xi'(u)}{u}\Bigr).
	\end{align*}
	Note that, since $\gamma_p>0$ for some $p\geq 3,$ it follows that
	\begin{align*}
	\Bigl(\frac{\xi'(s)}{s}\Bigr)'&=\frac{\xi''(s)s-\xi'(s)}{s^2}=\frac{1}{s^2}\sum_{p\geq 2}\bigl(p(p-1)-p\bigr)\gamma_p^2s^{p-1}>0
	\end{align*}
	and consequently, $\int_{0}^u\bigl(\beta^2\xi'(s)-Cs\bigr)ds<0.$ However, by making a choice of 
	$$
	\alpha(s)=\alpha_\beta(u)1_{[0,u)}(s)+\alpha_\beta(s)1_{[u,1]}(s)
	$$
	in \eqref{prop1:proof:eq1}, we obtain that 
	\begin{align*}
	\alpha_\beta (u)\int_{0}^u\bigl(\beta^2\xi'(s)-Cs\bigr)ds\geq 0.
	\end{align*}
	This gives a contradiction since $\alpha_\beta(u)>0.$ Thus, $u=0$. By (\ref{prop1:eq1}), $\mu_\beta\not = \delta_0$, so the second atom $v>0$ carries some weight and $0<m<1$.
\qed

\medskip
\noindent
\textbf{Proof of Proposition \ref{prop2}.}
We will verify that $\alpha_\beta(t)$ is the Parisi measure using the characterization in Proposition 2.1 in \cite{TalagrandSphh} and adapting a similar argument as the examples discussed after Proposition 2.2 in \cite{TalagrandSphh}. From the assumptions on $\xi$, the function
$$
\varphi(t) = \frac{1}{\beta \xi''(t)^{1/2}}
$$
is decreasing, concave, $\varphi(0)<1$ and $\varphi(1)>0.$ Therefore, the equation $\varphi(q)=1-q$ has the unique solution $q$ in $(0,1).$ Furthermore, 
$$
- \varphi'(t) = \frac{\xi'''(t)}{2\beta\xi''(t)^{3/2}}
\,\,\mbox{ is non-decreasing and $-\varphi'(q)<1$.} 
$$
Therefore, the function
		\begin{align*}
		\alpha(t)&=\left\{
		\begin{array}{ll}
		-\varphi'(t),&\mbox{if $t\in[0,q)$},\\
		1,&\mbox{if $t\in[q,1]$},
		\end{array}
		\right.
		\end{align*}
defines a cumulative distribution function on $[0,1].$ A direct computation gives that
	\begin{align*}
	\hat{\alpha}(t)&:=\int_t^1\!\alpha(s)\,ds=\left\{
	\begin{array}{ll}
	\varphi(t),&\mbox{if $t\in[0,q)$},\\
	1-t,&\mbox{if $t\in[q,1]$}.
	\end{array}
	\right.
	\end{align*}
Let us define
    \begin{align*}
	F(t)=\beta^2 \xi'(t)-\int_0^t\frac{ds}{\hat{\alpha}(s)^2}
	\,\,\mbox{ and }\,\,
	f(t)=\int_0^t\! F(s)\, ds.
	\end{align*}
Observe that, for $0\leq t\leq q,$
	$
	F(t)
	=\beta^2\xi'(t)-\beta^2(\xi'(t)-\xi'(0))=0,
	$
because of the assumption that $\gamma_1=0$. On the other hand, for $t\in(q,1],$
	\begin{align*}
	F(t)&=\beta^2\xi'(t)-\beta^2\xi'(q)-\int_{q}^t\frac{ds}{(1-s)^2}=\int_{q}^t\Bigl(\beta^2\xi''(s)-\frac{1}{(1-s)^2}\Bigr)ds<0,
	\end{align*}
because the integrand is negative on $(q,1]$ by the assumption of $\xi$ and the definition of $q$. As a result, $\sup_{t\in [0,1]}f(t) = 0$ and $\{t\mid f(t)=0\} = [0,q].$ Since $[0,q]$ is the support of the probability measure $\mu$ with the c.d.f. $\alpha$, Proposition 2.1 in \cite{TalagrandSphh} implies that $\mu$ is the Parisi measure.
\qed

\end{document}